\documentclass[12pt,a4paper]{article}       

\usepackage{xcolor}
 \usepackage[skins,theorems]{tcolorbox}
\tcbset{highlight math style={enhanced,
  colframe=red,colback=white,arc=0pt,boxrule=1pt}}
  \usepackage[bookmarksopen, bookmarksnumbered, bookmarksopenlevel=2]{hyperref}
 \usepackage[UKenglish]{babel}
 \usepackage[toc,page]{appendix}
 \usepackage{amsmath}
 \usepackage{amssymb}
 \usepackage{graphicx}
 \usepackage{hhline}

 \usepackage[bf]{caption}
\usepackage{cite}
\usepackage[vcentermath]{youngtab}
\usepackage{geometry}
\usepackage{slashed}
\usepackage{stackrel}
\usepackage{mathtools}
\usepackage{cancel} 
\usepackage{multirow}
\usepackage[margin=0pt,font=small,labelfont=normalfont,skip=22pt]{subcaption}

\usepackage{empheq}
\usepackage{arydshln}

\newcommand{\bqa}{\begin{eqnarray}}
\newcommand{\eqa}{\end{eqnarray}}

\def\del{\partial}

\newenvironment{eqn*}{\begin{equation*}\begin{aligned}}{\end{aligned}\end{equation*}\noindent}
\hypersetup{
    pdftitle={},
    pdfauthor={},
    pdfsubject={}
}
\numberwithin{equation}{section}
\numberwithin{table}{section}

\makeatletter

\definecolor{dgreen}{rgb}{0,0.45,0.2}
\definecolor{dblue}{rgb}{0,0.0,0.5}

\newcommand{\be}{\begin{equation}}
\newcommand{\ee}{\end{equation}}
\newcommand{\beq}{\begin{equation}}
\newcommand{\eeq}{\end{equation}}
\newcommand{\ba}{\begin{aligned}}
\newcommand{\ea}{\end{aligned}}

\newcommand{\bea}{\begin{eqnarray}}
\newcommand{\eea}{\end{eqnarray}}

\newcommand\bi{\begin{itemize}}
\newcommand\ei{\end{itemize}}

\def\unit{{1\kern-.65ex {\rm l}}}
\def\1{{1\kern-.65ex {\rm l}}}

\newcount\hour \newcount\minute
\hour=\time \divide \hour by 60
\minute=\time
\def\now{%
\ifnum \hour<13
  \ifnum \hour=0 \advance \hour by 12 \number\hour:\else \number\hour:\fi%
     \ifnum \minute<10 0\fi%
     \number\minute%
\ A.M.%
\else \advance \hour by -12 \number\hour:%
  \ifnum \minute<10 0\fi%
  \number\minute%
  \ P.M.%
\fi%
}

\makeatother

\begin{document}

\begin{titlepage}
\begin{center}
\rightline{\small }

\vskip 15 mm

{\large \bf
Chern-Simons Theory, Holography and Topological Strings} 
\vskip 11 mm

 Cumrun Vafa

\vskip 11 mm

{\it Jefferson Physical Laboratory, Harvard University, Cambridge, MA 02138, USA}

\end{center}
\vskip 17mm

\begin{abstract}
In this note we present a brief overview of connections between Chern-Simons theory and topological strings.  A prominent role in this link has been played by large N dualities and holography.  We demystify this by explaining why the Kahler form should be viewed as dual to the field strength associated with a 3-form gauge potential, sourced by Lagrangian D-branes.  We explain how this leads to the computation of topological string amplitudes in terms of topological vertex for toric Calabi-Yau threefolds.  Furthermore, applications of topological strings to a conceptual derivation of Skein relations for link invariants as well as some of its physical applications to black hole physics are also reviewed.

({\it Contribution to the AMS Bulletin in honor of James Simons})
\end{abstract}

\vfill
\end{titlepage}

\section*{James Simons}
James Simons was not only a first rate mathematician whose work left a deep impact in both mathematics and physics but also his visionary philanthropic activities together with Marilyn Simons has had transformative impact in support of science.
I have been fortunate to have also known Jim as a friend.   Moreover I have many fond memories of our discussions in connection with the establishment and growth of the Simons Center for Geometry and Physics in Stony Brook.

\section{Introduction}
In this contribution to the AMS Bulletin I highlight one aspect of the mathematical work of James Simons, together with Shiing-Shen
Chern \cite{Chern:1974ft}.  In particular they define a circle valued functional defined on manifolds with a gauge bundle and gauge connection. 
Chern-Simons term has played a key role in many diverse areas in theoretical physics.  In this note I focus  on its relation with topological (A-model) string theory also known as Gromov-Witten invariants (see e.g., \cite{mirror_symmetry_ams}) and holography (large N dualities) \cite{Gopakumar:1998ki}, and how this link leads to a solution of all topological string amplitudes on toric Calabi-Yau threefolds \cite{Aganagic:2003db}.
This relates Chern-Simons theory with various aspects of string theory, enumerative geometry, and low dimensional topological invariants. This link exemplifies a dual role of Chern-Simons theory as both a mathematical framework and a physical tool for exploring deep structures of physical theories.

The main aim of this note is to present a brief overview of this link.  One aspect which has not appeared in the literature before is that in this paper we present a simple explanation of the physical/mathematical mechanism underpinning the large N duality of Chern-Simons theory for $SU(N)$ gauge bundles.  Even though there has been heuristic explanations of this phenomenon, here we present an argument for it which is more in line with other examples of holography found in string theory, and in particular how a large gauge field flux generated by a charged object leads to a geometric transition and a dual holographic description.  Moreover we explain how this duality leads to a solution of (A-model) topologial strings on local toric Calabi-Yau threefolds.
We also review some other applications of topological strings including a conceptual explanation of the skein relations for links, as well as its relation to microstates of black holes.

\section{Chern-Simons Theory}
The Chern-Simons Functional $CS(A)$ for a gauge connecection $A$ in a $d=2n-1$ dimensional manifold $M$ is defined \cite{Chern:1974ft} by the condition that $d CS(A)\propto Tr F^n$.  Even though this terms is not gauge invariant, once integrated over the manifold it leads to a well-defined circle-valued function which depends on the gauge connection on $M$. 
The proportionality factor above is determined by making this correspond to the Chern class of the bundle. This functional leads to the  defintion of a physical action.
\subsection{The Chern-Simons Action}

The Chern-Simons theory is a gauge theory for three dimensional manifolds motivated by the Chern-Smions functional $CS(A)$ with the action defined by \cite{Witten:1988hf}:
\[
S_{\text{CS}} = \frac{k}{4\pi} \int_M \text{Tr} \left( A \wedge dA + \frac{2}{3} A \wedge A \wedge A \right),
\]
where \( A \) is the gauge field connection for a gauge group $G$, \( k \), an integer, is the coupling constant (or level), and \( M \) is a three-dimensional manifold, and $\text{Tr}$ denotes the group invariant trace. This action is topological, meaning that it depends only on the global properties of \( M \), not on its specific geometry, as the metric does not feature in the definition of the theory. The partition funciton of the theory $Z$ is defined by integrating over the infinite dimensional space of gauge in-equivalent connections on $M$:
$$Z(M,k,G)=\int DA\  {\rm{exp}} (i S_{CS})$$
Notice that the action is gauge invariant, as large gauge transformations can shift $S_{CS}$ by integer multiples of $2\pi k$.
The motivation for defining it for the case of $d=3$ is that only in this dimension the leading term is quadratic in the gauge field $A$, and this leads to a well-defined perturbation theory expansion on $Z$ in the weak coupling limit where $k\rightarrow \infty$.  Note that this defines invariants\footnote{Due to quantum effects, one actually needs a framing for the three manifold to make this well-defined.}for three dimensional manifolds which depends only on $k$ and the gauge group $G$.  This invariant is known as Witten-Reshetikhin-Turaev (WRT) invriant.

\subsection{Knot Invariants from Chern-Simons Theory}

Using Chern-Simons theory, in addition to invariants for 3-manifold, one can define invariants for knots and links. For example, consider a knot $K\subset M$ (i.e. an embedding of $S^1$ in $M$.  Then consider the expectation value of gauge invariant objects such as trace of gauge holonomies around $K$ in some representation $R$:
$$W^R_K={\rm Tr_R} \ {\rm Pexp} (i \int_M A).$$
The invariant for the knot is defined by inserting this in the path-integral over the gauge connection defining the Chern-Simons theory.  The result of this averaging $P_K$ defined by 
\[
P^R_K = \langle W^R_K \rangle ={1\over Z}\int DA\  W^R_K \cdot {\rm{exp}} (i S_{CS}),
\]
leads to a knot invariant for $K$ in $M$. This connection reveals how Chern-Simons theory bridges quantum field theory and the combinatorial world of knot theory. One can generalize this to links by replacing the $W^R_K$ by the product of the various knots $\prod_i K_i$ and representations $R_i$ where $K_i$ are the components of the links.

\section{Gromov-Witten Invariants and Topological Strings}

It seems the story of enumerative geometry, and in particular counting holomoprhic curves embedded in manifolds is a totally unrelated topic to Chern-Simons theory and invariants for links and three manifolds.  However, surprisingly they are deeply related.  This connection leads to a string theory perspective on Chern-Simons theory, as we will now review.

Gromov-Witten invariants (see \cite{mirror_symmetry_ams} for a review), involves the study of holomprhic maps from curves to symplectic target manifolds.  The case of most interest in physics is when the target manifold is a Calabi-Yau threefold.  Counting such holomorphic maps plays a central role in string compactifications as it leads to robust computations for physical amplitudes.  Roughly speaking one is interested in computing
$$
N_{g,d} = \text{Number of genus-} g \text{ holomorphic maps of degree } d.
$$
$$F(t,\lambda)=\sum_{g,d} \lambda^{2g-2} e^{-t\cdot d}N_{g,d}$$
where degree $d$ corresponds to a class $d\in H_2(CY,\mathbb{Z})$, i.e., the homology class of the image of the holomorphic curve, $t$ denotes the Kahler class, $\lambda$ the string coupling constant and $F$ the free energy of the theory where $Z={\rm exp}(-F)$ is the partition function of the topological string.  The formal dimension for the space of holomorphic maps from a Riemann surface to Calabi-Yau threefolds is zero, which means morally $N_{g,d}$ it should be counting isolated holomorphic maps.  However, they often comes in families, but even in this case there is a canonical way to calculate $N_{g,d}$ when they come in families. These lead to $N_{g,d}$ as being fractions rather than integer for such cases, but an integer can nevertheless be defined (the GV invariants) \cite{Gopakumar:1998ii,Gopakumar:1998jq} which is counting roughly holomorphically embedded objects as opposed to maps (see also \cite{Maulik:2016rip} and references therein for a mathematical exposition).

The connection between these invariants and topological string theory provides a powerful tool for studying string compactifications and their implications for 4-dimensional physics, where the space time is viewed as $R^4\times M$.

\subsection{Open Topological Strings}
One can extend the study of topological strings to the case where the Riemann surface has boundaries, by requiring that the boundaries of the Riemann surface map to specified Lagrangian cycles--Lagrangian `D-branes'--in the target manifold modeled after the upper half-plane ending on the real line.  It was shown in \cite{Witten:1992fb}
that if we consider a non-compact Calabi-Yau threefold given by $T^*L$ where $L$ is a Lagrangian submanifold, and if we consider $N$ Lagrangian D-branes wrapping $L$, the study of open Gromov-Witten theory gets related to $SU(N)$ Chern-Simons theory on $L$.  Namely, the perturbative expansion of $SU(N) $ Chern-Simons gauge theory inolves ``fat graphs'' or ``ribbon graphs'' (``tHooft diagrams) on $L$ which can be directly related to the corresponding worldsheet diagrams of the open Gromov-Witten theory.  The ribbon graphs arise as the gauge connection is a 1-form valued in the adjoint of the $SU(N)$, and the adjoint
can be viewed as a tensor product of fundemental (one boundary of the ribbon) with its conjugate (the other boundary of the ribbon).  The $AdA$ term in the CS action leads to parts of the graph with the topology of ribbons with topology $I\times R$ where $I$ is an interval,  and the $A^3$ term leads to the trivalent regions of the ribbon graph where three ribbons meet.  In this relation, the string coupling constant $\lambda$ (which counts the Euler class of the Riemann surface) gets mapped to $1/(k+N)$ where $k$ is the level of the CS theory.  As we will note below, this can also be derived from string field theory, which itself can be viewed as a Chern-Simons theory on the loop space.
\subsection{Topological A Model vs. Toplogical B Model}

Calabi-Yau manifolds were conjectured to come in mirror pairs $(M,M')$ \cite{Lerche:1989uy}, whose hodge numbers flip: $h^{p,q}(M)=h^{d-p,q}(M')$.  Evidence for this conjecture has accumulated over the years (see \cite{mirror_symmetry_ams}) and a general strategy for derivation of it has been proposed for geometric cases \cite{Strominger:1996it} and more generally in \cite{Hori:2000kt}.  This symmetry exchanges the Kahler structure and complex structure of Calabi-Yau pairs.  Moreover it maps the topological strings, which is called topological A-model, to a mirror topological strings, which is call the topological B-model.  The topological B-model captures the variations of the complex structure and has been studied intensively beginning with the work \cite{Bershadsky:1993cx}.  The open version of A-model, involving D-branes wrapping Lagrangian subspaces, gets mapped to open version of topological B-model involving D-branes wrapping holomoric subspaces in CY.
For more detail see \cite{mirror_symmetry_ams}.  An important aspect of the A and B models is that they do not affect one another.  In other words, the A-model is invariant under the complex deformations and the B-model is invariant under the Kahler structure changes.  This decoupling has implications for deriving the Skein relations for links, as we will review later in this paper.
\section{Large \( N \) Dualities and Topological Strings}

\subsection{SU(\( N \)) Chern-Simons Theory on \( S^3 \)}
We now specialize to the case where we have N-copies of the D-brane on $L=S^3$ and the CY is given by $T^*S^3$, and so considering open topological strings leads to the study of $SU(N)$ Chern-Simons theory on $S^3$.
Moreover we consider the limit in which $N\gg 1$.

In the large \( N \) limit, this theory rather surprisingly becomes equivalent (dual) to the topological string theory on the resolved conifold, i.e. local CY 3-fold given by $O(-1)\oplus O(-1)$ bundle over $P^1$ without any Lagrangian branes, where the volume $t$ of $P^1$ is identified with $t=N\lambda$.  In other words the Gromov-Witten theory of {\it open} strings involving Riemann surfaces with boundaries on conifold with $N$ Lagrangian branes on $S^3$, becomes equivalent to the {\it closed} Gromov-Witten theory involving holomorphic maps from Riemann surfaces without boundaries \cite{Gopakumar:1998ki}
\[
Z_{\text{CS}}(S^3) \longleftrightarrow Z_{\text{top}}\big[ O(-1)+O(-1)\rightarrow P^1\big].
\]

\begin{figure}[ht]
    \centering
    \includegraphics[width=0.9\textwidth]{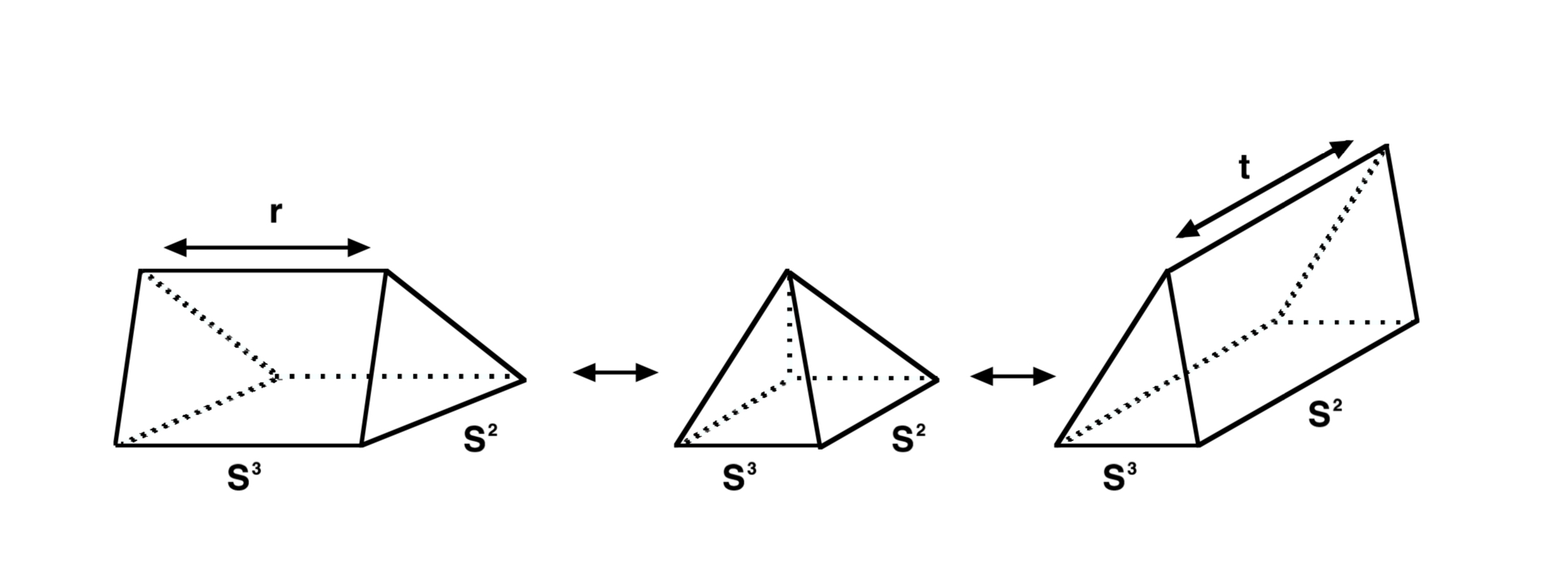}
    \caption{A large number $N$ of Lagrangian D-branes wrapped around $S^3$ leads to a geometric transition where the branes disappear and are replaced by the flux on the $S^2$ which links the $S^3$.  The Kahler class gets related to flux and leads to $N\lambda$ as the area of $S^2$, where $t=N\lambda$ is the string coupling.}
    \label{fig:transition}
\end{figure}
This has been checked by computing both sides independently and seeing that they agree.  The question is to better understand {\it why} there is this duality.  The key in understanding this, is the notion of holography in gravitational theories which induces geometric transitions which we discuss next.
\subsection{Geometric Transitions and Flux/Brane Duality: An Example of Holography}

The underpinning of the above duality is geometric transitions, where a contractible three-cycle with $N$ Lagrangian branes on it  shrinks and is replaced by a dual geometry without branes.
The key ingredient in this transition is that gravity (closed string modes) interact with gauge theory degrees of freedom (open string modes).  The imprint of the gauge degrees of freedom are the fluxes it gives rise to, and gravity responds to this by going through a geometric transition.  In particular the $P^1$ in the resolved conifold geometry, is already part of the geometry before transition, if we excise the loci of Lagrangian branes, i.e. the $S^3$.  Namely, the $P^1$ is nothing but the $S^2$ linking the $S^3$.  Indeed both $S^3$ and $S^2$ coexist there and at the geometric transition point where $S^3$ has shrunk we have a cone over $S^3\times S^2$ (see Figure \ref{fig:transition}).  Depending on which way we desingularize, by making the apex of the cone become $S^3$ or $S^2$ we go from the conifold geometry to the resolved conifold geometry.  This transition is topologically familiar in the context of Morse theory, as we change the level set of a Morse function passing through a crtical value.   

This type of transition provides a mechanism for understanding the dynamics of branes and fluxes in string theory and is essential in studies of gauge/gravity duality.  It provides an example of what is called `holography'.  In the present context it is the statement that in the theory without branes (a pure gravity theory) the boundary of the space encodes non-gravitational gauge degrees of freedom that encode what is going on inside.  In our example, this refers to the statement that the boundary of the resolved conifold includes an $S^3$ and the Chern-Simons theory on it, encodes the gravitational theory inside.

But how does this work mathematically, and what is the gravitational theory in the context of topological string?  How does $S^2$ get its size from the branes wrapping $S^3$. We now turn to answering these questions.

\section{Explanation of Holography Based on String Field Theory}

The relation between the description of the theory from the Riemann surface perspective (the worldsheet) vs. the target space description is captured by the string field theory.  In the case of open strings, the string field theory is of Chern-Simons type \cite{Witten:1985cc} given by
$$S=\frac{1}{\lambda}\int \frac{1}{2}\phi Q \phi +\frac{1}{3}\phi^3$$
where in the target space perspective for the conifold example we have the identification of the field:

$$\phi \leftrightarrow A\qquad Q\leftrightarrow  d$$

leading to the Chern-Simons action\footnote{Up to a renormalization of the CS level $k\rightarrow k+N$ and factors of $2\pi$ which we are not careful with here.}. 

For the closed string case the analog of the action is a bit more complicated but the quadradic term is rather simple:
$$S=\frac{1}{\lambda^2}\int{1\over 2} \Phi c_0^-Q \Phi+interactions $$
We now explain geometrically what corresponds to $c_0^-,Q, \Phi$ in the target space \cite{Bershadsky:1994sr}
for closed A-model topological strings.

For A-model topological string $\Phi$ corresponds to a $(1,1)$ form on the manifold, which could be viewed as the variation of the Kahler form from some fixed value $k_0$, i.e. $k=k_0+\Phi$. $Q$, just as in the open string case corresponds to the $d$ operator acting on forms. The operator $c_0^-$ is dual to the anti-ghost operator $b_0^-$, which is given by $d_c^\dagger=\partial^\dagger -\partial_c^\dagger$, being its dual implies that
$c_0^-={1\over d_c^\dagger}$.
We thus have the dictionary

$$\Phi \leftrightarrow C,\quad Q\leftrightarrow d,\quad b^-_0 \leftrightarrow d_c^{\dagger},\quad c_0^-\leftrightarrow \frac{1}{d_c^{\dagger}}.$$

However, the operator $c_0^-$ only makes sense if we delete the harmonic forms.  In particular $\Phi$ should have no component along those direction and can be written in the form
$$\Phi=d_c^\dagger C, \qquad k=k_0+d_c^\dagger C$$
In terms of this new field $C$ which is a 3-form the action becomes
$$S=\frac{1}{\lambda^2}\int \frac{1}{2}\Phi c_0^-Q \Phi+...=\frac{1}{\lambda^2}\int \frac{1}{2}d_c^{\dagger} C d C+...$$

Note that $C$ should be viewed as a higher form gauge field in the sense that $C$ is well-defined only up to a gauge transformation
$$C\rightarrow C+d_c^\dagger \epsilon ,$$
where $\epsilon$ is a 4-form gauge transformation parameter.

We are now ready to discuss the interaction of the open string field $C$ with the Lagrangian branes, which lead to the open string sector of the theory.
Consider a Lagrangian submanifold $L$ where we have wrapped a D-brane on it.  To see how this interacts with the closed string mode we should consider the disk amplitude where the boundary of the disk is restricted to lie on $L$ and we insert an operator $c_0^- \Phi$ at the center of the disk \cite{Sen:2024nfd}.  This insertion is nothing but the insertion of the 3-form field $C$ at the center and computation of this disk amplitude leads to a natural pairing of $C$ and $L$, i.e. $\int_L C$.  Taking into account the Euler character of the disk being 1, which fixes the power of $\lambda$, leads to the coupling
$$\frac{1}{\lambda}\int_{L}C$$
If we have $N$ D-branes wrapping $L$ this gives an extra factor of $N$ in the above, one for each pairing of the field C with the brane:
$$ \frac{N}{\lambda}\int_{L}C  $$

Putting this in the action we arrive at the form of the action including both the closed string sector, the open string gauge field and the interaction of the open and closed sectors:
$$S=\frac{1}{\lambda^2}\int_M \frac{1}{2}d_c^{\dagger} C \wedge d C+...+\frac{N}{\lambda}\int_{L}C+ \frac{1}{\lambda} \int_L \text{Tr} \left( \frac{1}{2}A \wedge dA + \frac{1}{3} A \wedge A \wedge A \right).$$

The coupling of the closed string field $C$ to the D-branes wrapping L, is of the same type as one has for electric charged objects sourcing Maxwell's fields.  The only difference being that everything is shifted upwards in dimension by 2 units.  To see this feature, note that if one varies the action with respect to $C$ the equation one gets is
$$d d_c^{\dagger} C=N\lambda\ \delta^3_{L}\rightarrow d (k-k_0)=dk=N\lambda\ \delta^3_{L}$$
From this we in particular learn that the Lagrangian D-branes source the Kahler form and the Kahler form is no longer closed unless we excise the region including the Lagrangian brane.
Moreover if we consider a 2-cycle linking $L$, we learn that
$$\int k=N\lambda$$
For example consider the case of $L=S^3$ we have discussed before with N D-branes wrapping $L$.  This means the 2-sphere linking with $S^3$ which originally had $\int_{S^2} k=0$, should now pick up an area 
$$\int_{S^2}k=N\lambda $$
In other words, to have a Kahler geometry description of the theory with N Lagrangian branes, where $k$ is closed, is equivalent to going to a small resolution of the conifold where $S^2$ has picked up a size.  This would get rid of the $S^3$ as a non-trivial cycle, and together with it, we have no more Lagrangian branes.  So we learn that the $SU(N)$ Chern-Simons theory on $S^3$ should be equivalent to closed topological strings on the resolution of conifold with Kahler class of ${\bf P}^1$ being $N \lambda$.

This phenomenon can also be seen in the mirror description of the theory which relates to the quantization of Kodaira-Spencer complex structure deformation \cite{Bershadsky:1993cx}.   In the B-model the role played by Kahler form is played by the holomorphic 3-form, which controls the complex structure.  In particular the closed string action leads to the Kodaira-Spencer theory of gravity \cite{Bershadsky:1993cx}.  Similarly the role of Lagrangian D-branes is played by holomorphic curves.
We have

$$\Phi \leftrightarrow B, \qquad Q\leftrightarrow {\overline \partial}, \qquad b_0^-\leftrightarrow \partial,\qquad  c_0^-=\frac{1}{\partial} $$
$$\Phi =\del B, \qquad \Omega=\Omega_0+{\overline \partial} B$$

where $B$ is a 2-form.  
$$S=\frac{1}{\lambda^2}\int\frac{1}{2} \Phi c_0^-Q \Phi+... =\frac{1}{\lambda^2}\int\frac{1}{2} \partial B {\overline \partial }B+...$$
Again the coupling between the D-branes represented by holomorphic curves ${\cal C}$ and the closed string field $B$ can be obtained, using the $c_0^-\Phi=B$ insertion on a disk ending on the holomorphic curve, leading to the pairing 
$${1\over \lambda}\int_{\cal C}B$$

If we have $N$ holomorphic branes wrapping a curve ${\cal C}$, the variation of the action with respect to $B$ will lead to the equation

$${\overline \partial}\partial B=N\lambda \  \delta_{\cal C} \rightarrow d \Omega= N\lambda \  \delta_{\cal C}$$
This implies that for any three-cycle linking ${\cal C}$ the holomorphic 3-form integrates to 
$$\int \Omega =N\lambda$$
which implies that we will need a geometric transition for this to correspond to a holomorphic 3-form.  In particular if we have the holomorphic D-branes wrapping the ${\bf P}^1$ in the resolved conifold geometry, the integral of $\Omega$ over $S^3$ linking with the $S^2$ is $N\lambda$ and this leads to a geometric transition where ${\bf P}^1$ shrinks and disappears and the $S^3$ grows, i.e. we are back to the conifold geometry of $T^*S^3$.

\section{Topological Vertex }
The holographic duality between Large N Chern-Simons theory and topological strings can be used to compute all genus topological string amplitudes for toric Calabi-Yau threefolds \cite{Aganagic:2003db}.  Let us first illustrate this with the example of the topological strings on the local CY given by the $O(-3)$ line bundle over ${\bf P}^2$.  

\begin{figure}[ht]
    \centering
    \includegraphics[width=0.7\textwidth]{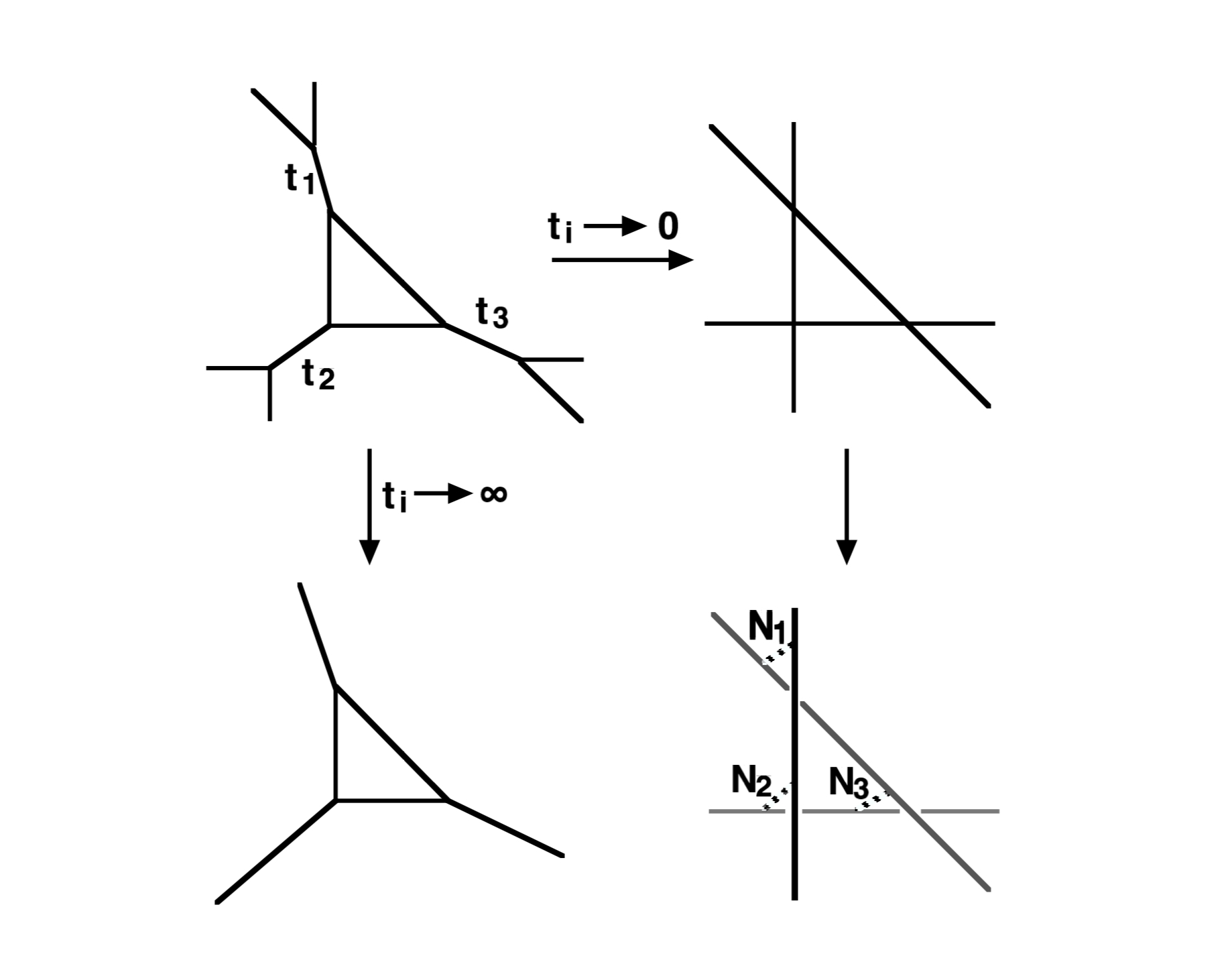}
    \caption{The topological Chern-Simons transitions, can be used to compute topological string amplitudes on ${\bf P}^2$ blown up at three points, which after three flops can be used to compute topological strings on ${\bf P}^2$ itself.}
    \label{fig:P2}
\end{figure}
To compute topological strings for this example, it turns out to be useful to consider a generalization of it, given by blowing up the ${\bf P}^2$ at three points.
If we flop the blown up curves to the fiber and take the corresponding Kahler parameters to infinity, we recover the ${\bf P}^2$ \ref{fig:P2}.  On the other hand the three times blown up geometry can be viewed as arising from three local $S^3\rightarrow S^2$
transitions as shown in Fig.\ref{fig:P2}.
So by large $N$ dualities this will have a reformulation in terms of the Chern-Simons theories $SU(N_1)\times SU(N_2)\times SU(N_3)$ where $t_i=N_i\lambda$ with some Wilson loop observables inserted along the legs of the toric diagram labeled by representations $R_1,R_2,R_3$
$$Z(t,t_i,\lambda)=\sum_{R_1,R_2,R_3} {\rm exp}(-t(d_1+d_2+d_3))\langle H^1_{R_1,R_2}H^2_{R_2,R_3}H^3_{R_3,R_1}\rangle $$
where $d_i$ are the degree of the representations $R_i$, $t$ is the Kahler class of ${\bf P}^2$ and $H^i(R_a,R_b)$ refer to the Hopf link correlation functions of the $SU(N_i)$ theory with representations $R_a,R_b$ on the two circles.

\begin{figure}[ht]
    \centering
    \includegraphics[width=0.7\textwidth]{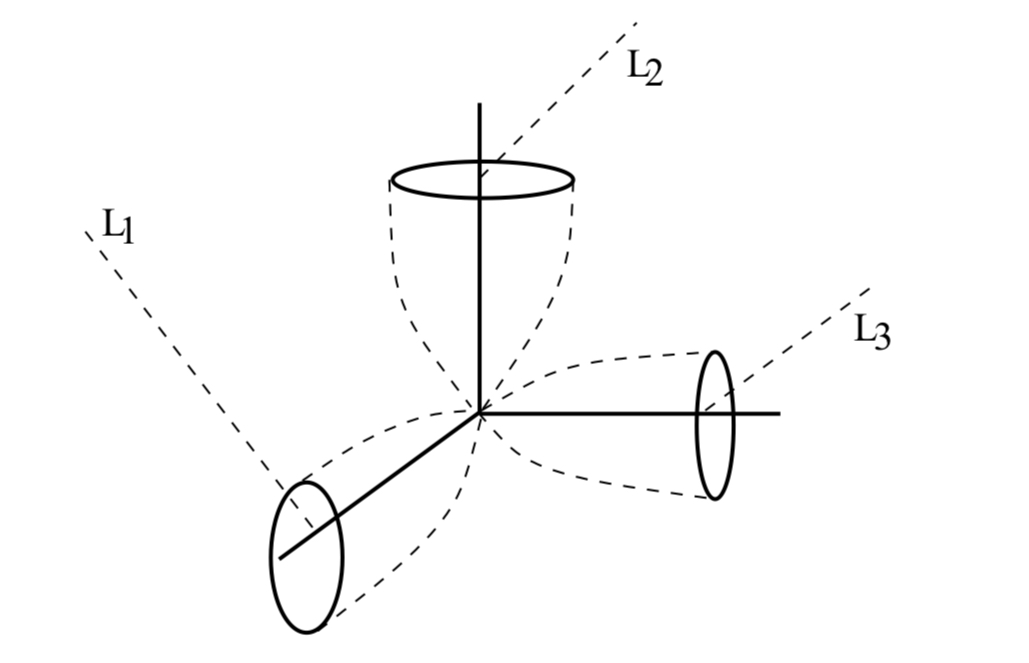}
    \caption{Topolotical vertex can be defined by holographic transitions from Chern-Simons theory.  It can be formulated in terms of three Lagrangian branes on ${\bf C}^3$ and the holomorphic curves ending on them can be labeled by Representations $R_i$ of $SU(N_i)$ for large $N_i$, leading to the vertex $C_{R_1,R_2,R_3}$.}
    \label{fig:top}
\end{figure}

This can be generalized to all toric CY.  Using these kinds of geometric transitions
 the computation of the topological A model can be reduced using localization and toric symmetries to a graphical recipe by putting together the local toric descriptions of ${\bf C}^3$, where each of the three legs of ${\bf C}^3$ has Lagrangian D-branes and the holomorphic curves ending on specific representation $(R_1,R_2,R_3)$ can be computed using the large N-duality between $SU(N)$ Chern-Simons theory on $S^3$.  This object which is called the topological vertex $C_{R_1,R_2,R_3}$ can be viewed as a computation on ${\bf C}^3$ with three Lagrangian branes $L_i$ supporting representations $R_i$ \ref{fig:top} on each leg
(for a mathematical exposition to topological vertex see \cite{Li:2004uf} and for a skein valued approach see \cite{Ekholm:2024lir}).

\section{Decoupling of A and B models and Skein Relations}
As already mentioned, one of the key features of topological strings is that there are two decoupled versions:  One is sensitive to the variations of Kahler structure and the other to the variations of complex structure.  In particular if you change the complex structure the computations of correlation functions involving the A-model do not change and similarly if we change the Kahler structure the amplitudes of the B-model do not change.
The fact these deformations do not affect each other has powerful consequences, and here we will note one such application, which has been recently used to derive the Skein relations for the HOMFLYPT polynomials \cite{Ekholm:2019yqp}.

The idea is to show that these skein relations are nothing but the statement of invariance of the amplitudes of the A-model with respect to variations of complex structure.  In particular one considers the A-model with $N$ D-branes wrapped on some Lagrangian branes and consider the contribution of string amplitudes to the Wilson loop observables.  As already noted, in the A-model these correspond to worldsheet diagrams ending on the Wilson lines in the D-branes.  These are called the worldsheet skein relations \cite{Ekholm:2024ceb} which lead to the Skein relations for the HOMFLYPT polynomials.
One considers variation of the complex structure and show that the following moves get induced (see \cite{Ekholm:2019yqp} for more detail):

\begin{figure}[ht]
    \centering
    \includegraphics[width=0.7\textwidth]{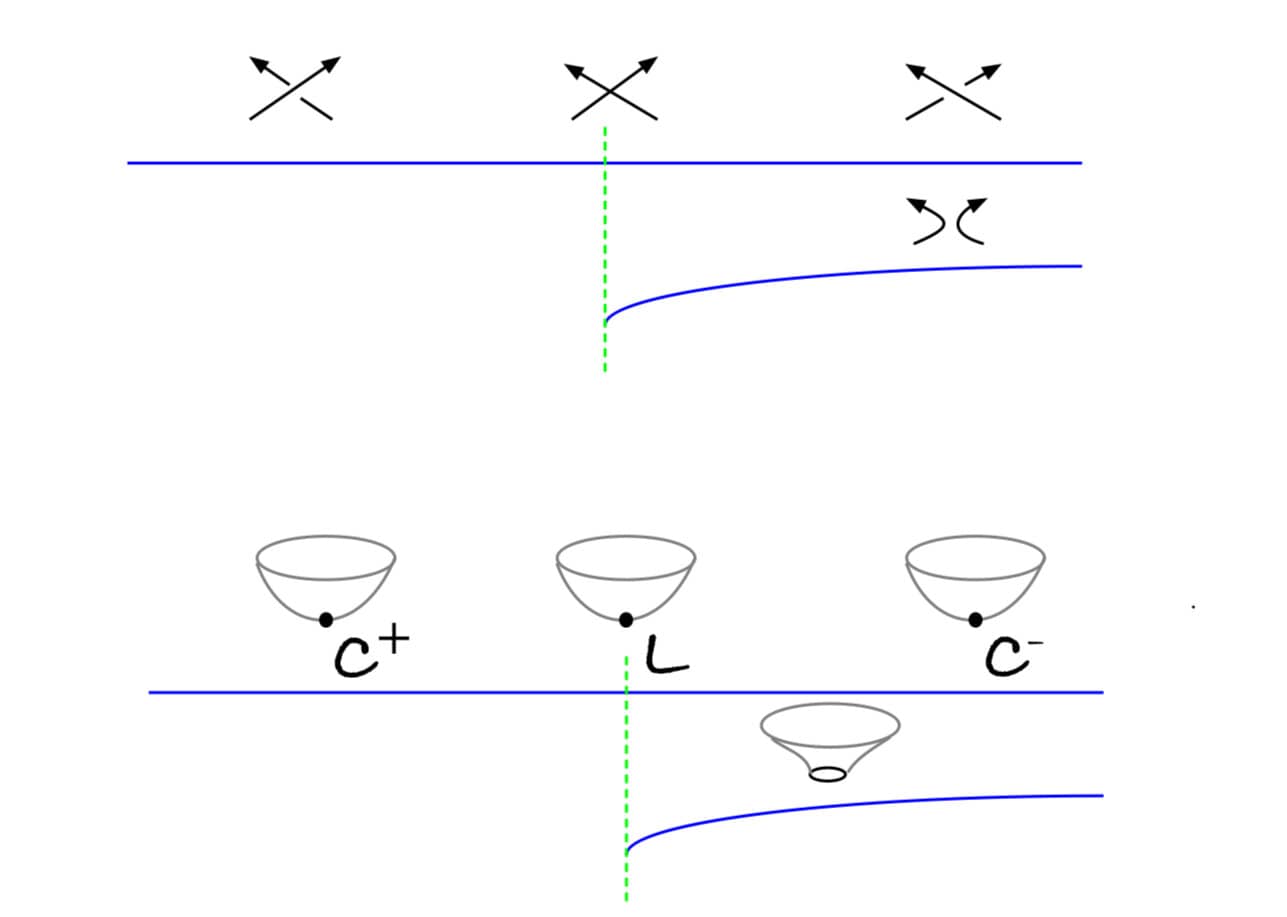}
    \caption{Worldsheet skein relations follow from recalling that changing complex structure should not change the amplitudes of the A-model. These figures represent change of complex structure from left to right.  The lines in the first figure represent the end points of the worldsheet diagrams, and it shows that under that complex deformations new worldsheet maps emerge after passing through a critical point (the dashed green line) where the image of the boundary of the worldsheet intersects. The second figure shows what happens if a Lagrangian D-brane crosses a closed holomorphic curve $C^+$ as one changes the complex structure.  Here it shows that new contributions emerge where the worldsheet can end on the Lagrangian brane as complex structure changes as well as a different Kahler class seen by the transformed holomorphic curve $C^-$ as its linking with $L$ has changed.}
    \label{fig:enter-label}
\end{figure}

These figures represent the move from left to right as we change the complex structure and follow what happens to the holomorphic curves in the CY as well as the possibility of getting new holomorphic curves.  The boundaries of the worldsheet correspond to contributions to Wilson line observables.  The equality of the A-model on the left with the A-model on the right incuded by the change of complex structure leads to the Skein relation for the HOMFLYPT  skein relations (for the fundamental representation), which in turn lead to complete computation of the HOMFLYPT invariants for knots and links.
An important ingredient in this relation is the second worldsheet skein move where the computation of the Kahler class on the closed holomorphic curve changes as we go from left to right, and the physical discussion above in terms of the Lagrangian branes giving rise to a shift in Kahler class due to the `Kahler flux' is a key ingredient.

\section{Physical Implications}
As we have already discussed Chern-Simons theory leads to a deeper understanding of topological strings, and in particular leading to a complete solution of topological strings for non-compact Calabi-Yau.  Topological strings on the other hand has had many applications in string theory \cite{Neitzke:2004ni}.  These include applications to partition function of gauge theories in various dimensions (see e.g.\cite{Iqbal:2012xm}).  This is done by geometrically engineering the relevant non-compact Calabi-Yau 3-fold which leads to the gauge theroy in the uncompactified spacetime.  The partition funciton of topological string on the corresponding Calabi-Yau, gets related to their partition function on $S^d$ or $S^{d-1}\times S^1$.

Topological strings can also be used also to count the microstates of BPS black holes. For example 5 dimensional BPS black holes can be represented by holomorphic curves embedded in the Calabi-Yau which lives in the compact 6 dimensional space, where the 11-dimensional M-theory is compactified on (note $11=6+5$).  On the other hand the counting of holomorphic curves on Calabi-Yau can be done via topological strings \cite{Gopakumar:1998ii,Gopakumar:1998jq,Katz:1999xq}.

Similarly for 4d BPS black holes one can use type IIB string theory compactifications from 10 to 4 on Calabi-Yau 3-fold.  In this case the BPS black holes can be viewed as D-branes wrapping special Lagrangian 3-cycles in the Calabi-Yau.  It turns out that in this case the leading count of the BPS black hole can be represented \cite{Ooguri:2004zv} by the square of the B-model topological strings on the Calabi-Yau 3-fold:
$$Z_{BH}=|Z^{top}|^2$$

\section{Conclusion}
I hope I have conveyed some of the rich applications that Chern-Simons theory has in string theory and how this leads to a deeper understanding of enumerative geometry and topological invariants for links.

\subsubsection*{Acknowledgments}

I would like to thank Tobias Ekholm for valuable discussions.

This work is supported in part by a grant from the Simons Foundation (602883,CV) and gifts from the DellaPietra Foundation.

\bibliographystyle{jhep}
\bibliography{sample}

\providecommand{\href}[2]{#2}\begingroup\raggedright\begin{thebibliography}{10}

\bibitem{Chern:1974ft}
S.-S.~Chern and J.~Simons, \emph{{Characteristic forms and geometric invariants}}, \href{https://doi.org/10.2307/1971013}{\emph{Annals Math.} {\bfseries 99} (1974) 48}.

\bibitem{mirror_symmetry_ams}
K.~Hori, S.~Katz, A.~Klemm, R.~Pandharipande, R.~Thomas, C.~Vafa et~al., \emph{Mirror Symmetry}, vol.~1 of \emph{Clay Mathematics Monographs}, American Mathematical Society and Clay Mathematics Institute, Providence, RI (2003).

\bibitem{Gopakumar:1998ki}
R.~Gopakumar and C.~Vafa, \emph{{On the gauge theory / geometry correspondence}}, \href{https://doi.org/10.4310/ATMP.1999.v3.n5.a5}{\emph{Adv. Theor. Math. Phys.} {\bfseries 3} (1999) 1415} [\href{https://arxiv.org/abs/hep-th/9811131}{{\ttfamily hep-th/9811131}}].

\bibitem{Aganagic:2003db}
M.~Aganagic, A.~Klemm, M.~Marino and C.~Vafa, \emph{{The Topological vertex}}, \href{https://doi.org/10.1007/s00220-004-1162-z}{\emph{Commun. Math. Phys.} {\bfseries 254} (2005) 425} [\href{https://arxiv.org/abs/hep-th/0305132}{{\ttfamily hep-th/0305132}}].

\bibitem{Witten:1988hf}
E.~Witten, \emph{{Quantum Field Theory and the Jones Polynomial}}, \href{https://doi.org/10.1007/BF01217730}{\emph{Commun. Math. Phys.} {\bfseries 121} (1989) 351}.

\bibitem{Gopakumar:1998ii}
R.~Gopakumar and C.~Vafa, \emph{{M theory and topological strings. 1.}},  \href{https://arxiv.org/abs/hep-th/9809187}{{\ttfamily hep-th/9809187}}.

\bibitem{Gopakumar:1998jq}
R.~Gopakumar and C.~Vafa, \emph{{M theory and topological strings. 2.}},  \href{https://arxiv.org/abs/hep-th/9812127}{{\ttfamily hep-th/9812127}}.

\bibitem{Maulik:2016rip}
D.~Maulik and Y.~Toda, \emph{{Gopakumar-Vafa invariants via vanishing cycles}},  \href{https://arxiv.org/abs/1610.07303}{{\ttfamily 1610.07303}}.

\bibitem{Witten:1992fb}
E.~Witten, \emph{{Chern-Simons gauge theory as a string theory}}, {\emph{Prog. Math.} {\bfseries 133} (1995) 637} [\href{https://arxiv.org/abs/hep-th/9207094}{{\ttfamily hep-th/9207094}}].

\bibitem{Lerche:1989uy}
W.~Lerche, C.~Vafa and N.P.~Warner, \emph{{Chiral Rings in N=2 Superconformal Theories}}, \href{https://doi.org/10.1016/0550-3213(89)90474-4}{\emph{Nucl. Phys. B} {\bfseries 324} (1989) 427}.

\bibitem{Strominger:1996it}
A.~Strominger, S.-T.~Yau and E.~Zaslow, \emph{{Mirror symmetry is T duality}}, \href{https://doi.org/10.1016/0550-3213(96)00434-8}{\emph{Nucl. Phys. B} {\bfseries 479} (1996) 243} [\href{https://arxiv.org/abs/hep-th/9606040}{{\ttfamily hep-th/9606040}}].

\bibitem{Hori:2000kt}
K.~Hori and C.~Vafa, \emph{{Mirror symmetry}},  \href{https://arxiv.org/abs/hep-th/0002222}{{\ttfamily hep-th/0002222}}.

\bibitem{Bershadsky:1993cx}
M.~Bershadsky, S.~Cecotti, H.~Ooguri and C.~Vafa, \emph{{Kodaira-Spencer theory of gravity and exact results for quantum string amplitudes}}, \href{https://doi.org/10.1007/BF02099774}{\emph{Commun. Math. Phys.} {\bfseries 165} (1994) 311} [\href{https://arxiv.org/abs/hep-th/9309140}{{\ttfamily hep-th/9309140}}].

\bibitem{Witten:1985cc}
E.~Witten, \emph{{Noncommutative Geometry and String Field Theory}}, \href{https://doi.org/10.1016/0550-3213(86)90155-0}{\emph{Nucl. Phys. B} {\bfseries 268} (1986) 253}.

\bibitem{Bershadsky:1994sr}
M.~Bershadsky and V.~Sadov, \emph{{Theory of Kahler gravity}}, \href{https://doi.org/10.1142/S0217751X96002157}{\emph{Int. J. Mod. Phys. A} {\bfseries 11} (1996) 4689} [\href{https://arxiv.org/abs/hep-th/9410011}{{\ttfamily hep-th/9410011}}].

\bibitem{Sen:2024nfd}
A.~Sen and B.~Zwiebach, \emph{{String Field Theory: A Review}},  \href{https://arxiv.org/abs/2405.19421}{{\ttfamily 2405.19421}}.

\bibitem{Li:2004uf}
J.~Li, C.-C.M.~Liu, K.~Liu and J.~Zhou, \emph{{A Mathematical theory of the topological vertex}}, \href{https://doi.org/10.2140/gt.2009.13.527}{\emph{Geom. Topol.} {\bfseries 13} (2009) 527} [\href{https://arxiv.org/abs/math/0408426}{{\ttfamily math/0408426}}].

\bibitem{Ekholm:2024lir}
T.~Ekholm, P.~Longhi and V.~Shende, \emph{{The skein valued mirror of the topological vertex}},  \href{https://arxiv.org/abs/2412.15454}{{\ttfamily 2412.15454}}.

\bibitem{Ekholm:2019yqp}
T.~Ekholm and V.~Shende, \emph{{Skeins on Branes}},  \href{https://arxiv.org/abs/1901.08027}{{\ttfamily 1901.08027}}.

\bibitem{Ekholm:2024ceb}
T.~Ekholm, P.~Longhi and L.~Nakamura, \emph{{The worldsheet skein D-module and basic curves on Lagrangian fillings of the Hopf link conormal}},  \href{https://arxiv.org/abs/2407.09836}{{\ttfamily 2407.09836}}.

\bibitem{Neitzke:2004ni}
A.~Neitzke and C.~Vafa, \emph{{Topological strings and their physical applications}},  \href{https://arxiv.org/abs/hep-th/0410178}{{\ttfamily hep-th/0410178}}.

\bibitem{Iqbal:2012xm}
A.~Iqbal and C.~Vafa, \emph{{BPS Degeneracies and Superconformal Index in Diverse Dimensions}}, \href{https://doi.org/10.1103/PhysRevD.90.105031}{\emph{Phys. Rev. D} {\bfseries 90} (2014) 105031} [\href{https://arxiv.org/abs/1210.3605}{{\ttfamily 1210.3605}}].

\bibitem{Katz:1999xq}
S.H.~Katz, A.~Klemm and C.~Vafa, \emph{{M theory, topological strings and spinning black holes}}, \href{https://doi.org/10.4310/ATMP.1999.v3.n5.a6}{\emph{Adv. Theor. Math. Phys.} {\bfseries 3} (1999) 1445} [\href{https://arxiv.org/abs/hep-th/9910181}{{\ttfamily hep-th/9910181}}].

\bibitem{Ooguri:2004zv}
H.~Ooguri, A.~Strominger and C.~Vafa, \emph{{Black hole attractors and the topological string}}, \href{https://doi.org/10.1103/PhysRevD.70.106007}{\emph{Phys. Rev. D} {\bfseries 70} (2004) 106007} [\href{https://arxiv.org/abs/hep-th/0405146}{{\ttfamily hep-th/0405146}}].

\end{thebibliography}\endgroup

\end{document}